\documentclass[a4paper, 10pt, leqno,final]{article}

\usepackage{eucal}	
\usepackage{amsmath}
\usepackage{amsfonts}
\usepackage{amsthm}
\usepackage{amssymb}
\usepackage{graphicx}
\usepackage{graphics}
\usepackage{bm}
\usepackage[perpage]{footmisc}

\usepackage{dsfont}
\usepackage{color}
\usepackage{enumerate}
\usepackage[font=footnotesize]{caption}
\usepackage[all]{xy}
\usepackage{xypic}
\allowdisplaybreaks
\usepackage{enumerate}		

\usepackage{mathtools}		
\usepackage{mathrsfs}		
\usepackage{eucal}			
\usepackage{braket}	
\usepackage[colorlinks=true, linkcolor=black, citecolor=black,
urlcolor=blue]{hyperref}		
\mathtoolsset{showonlyrefs}
\usepackage{mparhack}		
\usepackage{relsize}			
\usepackage{a4wide}		
\usepackage{booktabs}		
\usepackage{multirow}		
\usepackage{caption}		
\usepackage{rotating}
\usepackage{varioref}		
\usepackage{footmisc}

\numberwithin{equation}{section}

\theoremstyle{plain}
\newtheorem{thm}{Theorem}[section]
\newtheorem*{mainthm}{\textsc{Theorem A}}

\newtheorem{prop}[thm]{Proposition}
\theoremstyle{definition}

\newtheorem{dfn}[thm]{Definition}

\theoremstyle{remark}
\newtheorem{rmk}[thm]{Remark}
		\newtheorem{rem}[thm]{Remark}		
				\newtheorem{note}[thm]{Notation}	

		\renewcommand{\leq}{\leqslant}
\renewcommand{\geq}{\geqslant}
\renewcommand{\=}{\coloneqq}			

\newcommand{\Lag}{\mathrm{Lagr}}
\newcommand{\eul}{\mathrm{Eul}}
\newcommand{\N}{\mathds{N}}

\newcommand{\Z}{\mathds{Z}}
\newcommand{\Q}{\mathbb{Q}}

\newcommand{\R}{\mathds{R}}

\newcommand{\SO}{\mathrm{SO}}

\renewcommand{\O}{\mathrm{O}}

\newcommand{\sgn}{\mathrm{sgn}}

\DeclareMathOperator{\spfl}{sf}

\DeclareMathOperator{\diag}{diag}		

\newcommand{\cc}{\mathrm{CC}}

\newcommand{\sbc}{\mathrm{SBC}}
\newcommand{\csbc}{\mathrm{CSBC}}

\renewcommand{\=}{\coloneqq}			

\usepackage{bm}

\title{Bifurcations of balanced configurations\\ for the Newtonian $n$-body problem in $\R^4$}
\author{Luca Asselle, Marco Fenucci, Alessandro Portaluri}
\date{\today}

\date{\today}
\begin{document}
 \maketitle

\begin{center}
\large \textit{Dedicated to Claude Viterbo on the occasion of his sixtieth birthday}
\end{center}

\vspace{5mm}

\begin{abstract}
For the gravitational $n$-body problem, the simplest motions are provided by those rigid motions in which each body moves along a Keplerian orbit and the shape of the system is a constant (up to rotations and scalings) configuration featuring suitable properties. While in dimension $d \le 3$ the configuration must be \textit{central}, in dimension $d \ge 4$ new possibilities arise due to the complexity of the orthogonal group, and indeed there is a wider class of $S$-\textit{balanced configurations}, containing central ones, which yield simple solutions of the $n$-body problem. Starting from the recent results in \cite{AP:20}, we study the existence of continua of bifurcations branching from a trivial branch of collinear $S$-balanced configurations and provide an estimate from below on the number of bifurcation instants. In the last part of the paper, by using the continuation method, we explicitly display the bifurcation branches in the case of the three body problem for different choices of the masses.

\vspace{2mm}

\end{abstract}

%


\section{Introduction}\label{sec:introduction}

 The Newtonian $n$-body problem concerns the motion of $n$ point particles with masses $m_j \in \R^+$ and positions $q_j \in \R^d$, where $j=1, \ldots, n$ and $ d \ge 2$, interacting each other according to Newton's law of inverses squares. The particles thus move according to Newton's equations of motion, which in this case read 
 \begin{equation}\label{eq:motions}
 m_j \ddot q_j = \dfrac{\partial U}{\partial q_j} \quad \textrm{ where } \quad U(q_1, \ldots, q_n)\=\sum_{i<j}\dfrac{m_i m_j}{|q_i-q_j|}. 	
 \end{equation}
Letting $M$ be the $(nd \times nd)$-diagonal {\sc mass matrix\/} defined by 
$$M\= \diag(\underbrace{m_1, \ldots, m_1}_{d{\text{-times}}}, \ldots, \underbrace{m_n, \ldots, m_n}_{d{\text{-times}}})$$ the equations of motion can be equivalently written as
\begin{equation}\label{eq:Newton-intro}
	\ddot q= M^{-1}\nabla U(q).
\end{equation}
As the center of mass has an inertial motion, we can fix it without loss of generality at the origin.
Among all possible configurations of the system, a crucial role  is played by the so-called {\sc central configurations\/} ($\cc$  for short), namely by those configurations in which $M^{-1}\nabla U(q)$ is parallel to $q$: 
\begin{equation}\label{eq:cc-eq-2}
	M^{-1}\nabla U(q) + \lambda q=0.
\end{equation}   
In other words, the acceleration vector of each particle is pointing towards the origin with magnitude proportional to the distance to the origin. As a straightforward consequence of the homogeneity of the potential we obtain that the proportionality constant $\lambda$ is actually equal to $-U(q)/\langle Mq, q\rangle$. 

 Equation~\eqref{eq:cc-eq-2} is a non-linear algebraic equation which is almost impossible to solve explicitly, and despite substantial progresses (starting from the work of - among others - Smale, Conley, Albouy, Chenciner, McCord, Moeckel, Pacella) have been made in the last decades, many basic questions about $\cc$ still remain unsolved.
Nevertheless, there are several reasons why $\cc$ are of interest in the $n$-body problem and more generally in Celestial Mechanics: 
\vspace{-2mm}
 \begin{itemize}
\item[-] Every $\cc$ defines a {\sc homothetic solution\/} of~\eqref{eq:Newton-intro}, namely a solution which preserves its shape for all time while receding from or collapsing into the center of mass.
\vspace{-2mm}
\item[-] Planar $\cc$ give rise to a family of periodic motions of~\eqref{eq:Newton-intro}, the so-called {\sc relative equilibria\/}, in which the configuration rigidly rotates at a constant angular speed about the center of mass. 
More generally, any such $\cc$  gives rise to a family of  {\sc homographic  solutions\/} of~\eqref{eq:Newton-intro} in which each particle traverses an elliptical orbit  with eccentricity $ e \in (0,1)$.
\vspace{-2mm}
\item[-] $\cc$ control the qualitative behavior of total colliding solutions (and completely parabolic motions) of the $n$-body problem.
 \end{itemize}
 \vspace{-2mm}
For the $n$-body problem in $\R^d$, $d \le 3$, configurations which are not central cannot produce homographic motions at all. If we instead allow dimensions $d \ge 4$, then there is a wider class of so-called ``$S$-balanced configurations'' which produces relative equilibria of the $n$-body problem. These new high dimensional phenomena were first observed by Albouy and Chenciner in \cite{AC98} (cfr. also \cite{Moe14}) and are due to the higher complexity of the orthogonal group, which allows, for example, to rotate in two mutually orthogonal planes with different angular velocities, thus leading to new ways of balancing the gravitational forces with centrifugal forces. We shall notice that, in contrast with the case $d=2$, the resulting relative equilibria will be periodic in time only if the angular velocities are rationally dependent, and quasi-periodic otherwise. 
 
We now define $S$-balanced configurations rigorously in the case $d=4$, which we will focus for the rest of the paper onto. 
Thus, fix a positive real number $s> 1$ and consider the $(4\times 4)$-diagonal matrix 
$$S= \diag (s,s,1,1).$$
Any solution of 
\begin{equation}\label{eq:intro-s-balanced-conf}
M^{-1}\nabla U(q) + \lambda \widehat S q =0, \quad \widehat S := \text{diag}(S,...,S),
\end{equation}
is called an $S$-{\sc balanced configuration}, $\sbc$ for short, and gives rise to a uniformly rotating relative equilibrium solution of Equation~\eqref{eq:Newton-intro} in $\R^4$. Clearly, for $s=1$, we obtain again \eqref{eq:cc-eq-2}. Also, 
Equation~\eqref{eq:intro-s-balanced-conf} is, for every $s>1$, invariant under the non-free (actually, not even locally free) diagonal $S^1\times S^1$-action given by rotations in the $\R^2\times \{0\}$ and $\{0\}\times \R^2$ planes. In particular, solutions always come in families, namely in $S^1$-families if they are contained in one of the two planes above and in $S^1\times S^1$-families otherwise.

\begin{rmk}
For $n=3$, there is a big class of planar non-equilateral and non-collinear isosceles triangles which are $\sbc$ but not $\cc$; for further details we refer to \cite{Moe14}. From a physical viewpoint, the larger $s$ is the faster the bodies contained in the plane $\R^2 \times \{0\}$ rotates. Such a rough physical interpretation is paradigmatic of a deep stability issue which we are currently investigating. \qed
\end{rmk} 

In the study of  Equation~\eqref{eq:intro-s-balanced-conf} it is quite natural to interpret $s$ as a bifurcation parameter, and hence to try to understand if: 
\vspace{-2mm}
\begin{enumerate}[i)] 
\item there exist configurations which are $\sbc$ (possibly collinear) for every choice of the parameter $s$ (in other words, whether or not there are trivial branches of solutions), and 
\vspace{-2mm}
\item how many (if any) bifurcation points one has along such trivial branches.  
\end{enumerate}

\vspace{-2mm}

As far as Question i) is concerned, we readily see that collinear $\cc$ in the plane $\{0\}\times \R^2$ are solutions 
of \eqref{eq:intro-s-balanced-conf} independently of $s>1$ and hence define trivial branches of solutions $(\widehat{q}_s)_{s>1}$. 
Using the variational characterization of $\sbc$ (for more details see \cite{AP:20} or the following section) and the fact that along the trivial branches 
the Morse index jumps at precisely characterized values of the parameter $s$, we will provide the following answer to Question ii). For a more precise statement we refer to Theorem~\ref{main:sbalanced-bifurcation}.

\begin{mainthm}
For $s_1$ sufficiently close to $1$ and $s_2$ large enough, there are at least $n!$ bifurcation instants from the trivial branches 
of solutions $(\widehat q_s)_{s \in [s_1,s_2]}$ corresponding to collinear $\cc$ in the $\{0\}\times \R^2$-plane.  
\end{mainthm}

We shall notice already at this point that the non-trivial branches emanating from the trivial ones are genuine $\sbc$ (that is, not $\cc$). This will be clear from the construction, anyway this also follows from the fact that collinear 
$\cc$ are isolated as central configurations in virtue of Moeckel's $45^\circ$-theorem. 

As already mentioned, the main idea behind the result above is that, for variational problems in finite dimension (and, under suitable condition,  also in infinite dimension), bifurcation instants along some trivial branch are detected by the jump of the Morse index as soon as the trivial branch is degenerate only in finitely many points. From a technical viewpoint, one difficulty to overcome is to rule out the degeneracy due to the $S^1\times S^1$-symmetry of Equation~\eqref{eq:intro-s-balanced-conf}. This will be done by means of 
a reduction argument (see Section~\ref{sec:preliminaries}).

In case $n=3$, we will use numerical methods to provide a rather complete description of the non-trivial branches bifurcating from the trivial branches of collinear $\cc$. Already in such an easy case, 
we observe some very interesting and rather unexpected phenomena: besides a strong dependence on the choice of the masses (which we recall is not the case for $\cc$, as for any choice of the masses one has precisely 4 $\cc$ up to symmetry, namely the three Euler configurations, which are saddle points of $U$, and Lagrange's equilateral triangle, which is a global minimum of $U$), we e.g. observe the presence 
of connections between Lagrange's equilateral triangle and (some of the) Euler configurations through paths of $\sbc$ which are for any $s$ local minima of $U$, as well as of turning points along some of the non-trivial branches at which the Morse index jumps but from which no secondary branches originate. This suggests that for larger values of $n$ extremely interesting new phenomena might occur. We plan to study these aspects further in future work. 

We shall also mention that other trivial branches can be constructed from planar $\cc$ in the plane $\{0\}\times \R^2$. Since the Morse index jumps also along such branches, we should be able to find other bifurcations instants. However, the problem is here more complicated since the degeneracy due to the symmetry cannot be overcome by reduction, and hence a generalization 
of the abstract bifurcation result (see Theorem~\ref{thm:main-abstract}) to an equivariant setting is needed. We plan to address this issue in future work. 

We end up this introduction with a brief summary of the content of this paper: In Section~\ref{sec:preliminaries} we define 
$\sbc$ and discuss their basic properties. In Section~\ref{sec:sf} we briefly recall the definition of the spectral flow in a finite 
dimensional setting. In Section~\ref{sec:4} we prove an abstract bifurcation result from the trivial branch of a one parameter 
$\mathscr C^2$-family of functions on a finite-dimensional manifold and then apply it to the study of bifurcations of $\sbc$.
Finally, in Section~\ref{sec:5} we use numerical computations to study the non-trivial branches bifurcating
from a trivial branch of collinear $\cc$ in the case $n=3$.

\vspace{3mm}

\noindent \textbf{Acknowledgments.} Luca Asselle is partially supported by the DFG-grant 380257369 ``Morse theoretical methods in Hamiltonian dynamics''. Marco Fenucci is partially supported by the MSCA-ITN Stardust-R, Grant Agreement n. 813644 under the H2020 research and innovation program.


\section{S-balanced configurations in the n-body problem} \label{sec:preliminaries}

In this section we recall the definition of $S$-balanced configurations and their basic properties, referring to \cite{AP:20} and \cite{Moe14} for the details. For $n\geq 2$, we consider $n$ point-masses 
$m_1,...,m_n>0$, whose positions are denoted by $q_1,...,q_n\in \R^4$ respectively, and which are supposed to interact with each other according to Newton's law of inverse squares. 
Setting the {\sc mass matrix} $M$ to be the diagonal $(4n\times 4n)$-matrix 
$$M:= \text{diag}\, (m_1 I_4,..., m_n I_4),$$
where $I_4$ is the $4$-dimensional identity matrix, we readily see that the equations of motion read 
\begin{equation}
M \ddot q = \nabla U(q)
\label{NE}
\end{equation}
where $q=(q_1,...,q_n)\in \R^{4n}$ is the {\sc configuration vector} of the $n$ point-masses, $\nabla$ denotes the gradient in $\R^{4n}$, and $U$ is the {\sc Newtonian (gravitational) potential}
$$U(q) := \sum_{i<j} \frac{m_im_j}{|q_i-q_j|}.$$ 
The invariance of \eqref{NE} under translations implies in virtue of Noether's theorem that the center of mass 
$$\overline q := \sum_{i=1}^n m_i q_i$$
has an inertial motion, and hence it is not restrictive to fix it at the origin. Therefore, we can without loss of generality suppose that $U$ is defined over the space of {\sc collision free configurations with center 
of mass at the origin} 
$$ \widehat{\mathbb X} := \Big \{ q=(q_1,...,q_n) \in \R^{4n} \ \Big |\ \overline q =0, \ q_i\neq q_j \ \forall i\neq j \Big \}.$$
The set over which $U$ is not defined, namely
$$\Delta := \{ q \in \R^{4n}\ |\  \overline q=0\} \setminus \widehat{\mathbb X},$$ 
is called the {\sc collision set}. As it is nowadays well-known, Equation~\eqref{NE} is extremely hard to solve, and indeed a complete solution is possible only for $n=2$. Therefore, instead of trying to 
solve \eqref{NE} explicitly, one can try to look for (simple) solutions with prescribed behavior. 

The simplest possible solutions of \eqref{NE} one can think of are those given by \textit{homographic motions}, i.e. rigid motions in which the configuration
of the bodies remains constant (up to rotations and scalings) in time. If one makes such an Ansatz in dimension 2 or 3, then one finds that
the configuration of the bodies must be {\sc central}, namely a solution of 
\begin{equation}
M^{-1} \nabla U(q) + \lambda q =0.
\label{CC}
\end{equation}
In other words, any solution of \eqref{CC} gives, for suitable choice of the initial momentum, rise to a homographic motion. As it turns out, in this case, each body must then move along a Keplerian orbit.
In the particular case of the zero angular momentum Keplerian orbit, we retrieve the so called \textit{homothetic motions} in which all masses collapse simultaneously at the origin or recede 
from total collision. In case of the eccentricity zero Keplerian orbit instead, we retrieve the so called \textit{relative equilibria}, in which the configuration of the bodies rigidly rotates around the origin 
at uniform speed while keeping its size constant. In dimensions 2 and 3, there are no other possible homographic motions. In dimension 4 instead new possibilities arise due to the higher 
complexity of the orthogonal group $\O(4)$: Indeed, in $\R^4$ it is possible to rotate simultaneously in two mutually orthogonal planes with different angular velocities. This produces a new
balance between gravitational and centrifugal forces, thus yielding new periodic or quasi-periodic motions. Thus, for $d=4$ there is a wider class of configurations, the so called $S$-{\sc balanced configurations}, 
which contains central ones and provides simple solutions to \eqref{NE}. 

More precisely, fix a positive real number $s> 1$ and consider the $(4\times 4)$-diagonal matrix 
$$S= \diag (s,s,1,1).$$
Any solution of 
\begin{equation}
M^{-1}\nabla U(q) + \lambda \widehat S q =0, \quad \widehat S := \text{diag}(S,...,S),
\label{BC}
\end{equation}
is called an $S$-balanced configuration. Clearly, for $s=1$ we obtain \eqref{CC}. It is easily seen that any solution of \eqref{BC} yields a relative equilibrium solution of \eqref{NE}, 
$$q(t) := \left (\begin{matrix} e^{i\sqrt{s} t} & 0 \\ 0 & e^{it}\end{matrix} \right )\cdot q,$$
which will then be a periodic solution if $s\in \Q$ and a quasi-periodic solution otherwise. Notice also that \eqref{BC} is invariant under the (diagonal) $S^1\times S^1$-action on $\widehat{\mathbb X}$ 
given by rotations in the $\R^2\times \{0\}$ and $\{0\}\times \R^2$ planes, whereas \eqref{CC} is $\SO(4)$-invariant. Both actions are not free (actually, not even locally free). 

Taking the scalar product of \eqref{BC} with $q$ 
and using Euler's theorem, we see that the constant $\lambda$ appearing in \eqref{BC} must be equal to 
$$\frac{U(q)}{I_S(q)},$$
where $I_S(q) := \langle \widehat{S} M q,q\rangle$ is the $S$-{\sc weighted moment of inertia}, and as a direct consequence of the invariance under scalings of \eqref{BC}, we see that we can always 
normalize an $S$-balanced configuration to satisfy $I_S(q)=1$. It is therefore natural to introduce the {\sc collision free configuration sphere}
$$\widehat{\mathbb S} := \Big \{q\in \widehat{\mathbb X} \ \Big |\ I_S(q)=1\Big \},$$
and to consider only {\sc normalized} $S$-{\sc balanced configurations}, i.e. solutions of \eqref{BC} which are contained in $\widehat{\mathbb S}$. We notice that Equation \eqref{BC} on $\widehat{\mathbb S}$ 
reads 
\begin{equation}
M^{-1} \nabla U(q) + U(q) \widehat S q =0.
\label{BCnormalized}
\end{equation}
To simplify the notation, we will hereafter refer to solutions of \eqref{BCnormalized} simply as $S$-balanced configurations. In other words, all $S$-balanced configurations will be hereafter assumed to be normalized. 
We will also use the shorthand notation $\sbc$ instead of $S$-balanced configuration.

\begin{rmk}
The interest on $S$-balanced and central configurations goes far beyond the fact that they yield simple solutions of \eqref{NE}. Indeed, their properties turn out to be useful to understand 
the qualitative behavior of many other classes of solutions to \eqref{NE}, as e.g. colliding solutions. For more details we refer to \cite{AP:20}. \qed
\end{rmk}

\begin{rmk}
$S$-balanced configurations have been introduced in the late nineties by Albouy and Chenciner, see \cite{AC98}. There, and also in \cite{Moe14}, the matrix $S$ is supposed to be minus the square
of a skew-symmetric matrix. As it turns out, our definition of $S$-balanced configurations is completely equivalent to that in \cite{AC98}. Indeed, after replacing the standard basis of $\R^4$ 
with a (orthonormal) spectral basis of $S$, we can suppose $S$ to be in diagonal form. Also, the invariance under scalings of the problem implies that we can suppose $S$ to be of the form considered above. \qed
\end{rmk}

A key feature of $\sbc$ is that they admit a variational characterization: Indeed, a configuration vector $q\in \widehat{\mathbb S}$ is a $\sbc$ if and only if it is a critical point of the restriction of $U$ 
to $\widehat{\mathbb S}$, which with slight abuse of notation will be hereafter denoted also with $U$. 

The Hessian of $ U : \widehat{\mathbb S} \to \R$ at a critical point $q$ is the quadratic form on $T_q\widehat{\mathbb S}$ represented, with respect to the mass-scalar product $\langle M\cdot,\cdot\rangle$,  by the $(4n\times 4n)$-matrix 
$$H(q) = M^{-1}D^2 U(q) + U(q) \widehat S.$$
 A straightforward computation shows that $(i,j)$-th entry of $D^2U(q)$ is given by 
 $$D_{ij} = \frac{m_im_j}{r_{ij}^3} \big ( I_4 - 3 u_{ij}u_{ij}^T \big ), \quad \text{for}\ i \neq j, \qquad D_{ii} = -\sum_{j\neq i} D_{ij},$$
where as usual one sets
 $$r_{ij} := |q_i-q_j|, \quad u_{ij} := \frac{q_i-q_j}{|q_i-q_j|}.$$
 As we already observed, Equation~\eqref{BCnormalized} is $S^1\times S^1$-invariant, and hence $\sbc$ always appear in $S^1\times S^1$-families. In particular, the Hessian $H(q)$ is always degenerate 
 as a quadratic form. Since we do not want to work in a setting where a group of symmetries is acting, we proceed as 
 follows using what in \cite{AP:20} we called the \textit{reduction to (H1) argument}: the submanifold 
 $$\mathcal P_s:= \Big \{q \in \widehat{\mathbb S}\ \Big |\ q_k \in \{0\}\times \R^2 \times \{0\},\ \forall k=1,...,n\Big \}\subset \widehat{\mathbb S}$$
 of planar configurations in the plane $\{0\}\times \R^2 \times \{0\}$ is invariant under the gradient flow of $U$ (actually, any submanifold of planar configurations in some coordinate plane does). 
 Ignoring all vanishing components, thus identifying $\{0\}\times \R^2 \times \{0\}$ with $\R^2$, we see that \eqref{BCnormalized} on $\mathcal P_s$ reads 
 \begin{equation}
 M^{-1}\nabla U(q) + U(q) \cdot \text{diag}\, \left ( \left ( \begin{matrix} s & 0 \\ 0 & 1 \end{matrix}\right ), ..., \left ( \begin{matrix} s & 0 \\ 0 & 1 \end{matrix}\right )\right ) q =0,
 \label{BCplane}
 \end{equation}
 where with slight abuse of notation we denote the mass-matrix on $\mathcal P_s$ and the restriction of $U$ to $\mathcal P_s$ again with $M$ and $U$ respectively. The main advantage of such a reduction argument is that \eqref{BCplane} is no longer invariant under the $S^1\times S^1$-action, but rather only under the action of the discrete group $\Z_2\times \Z_2$ given by reflections along the main axes. 
In particular, solutions of \eqref{BCplane}, which are nothing else but planar $\sbc$ contained in the plane $\{0\}\times \R^2\times \{0\}\subset \R^4$, do not a priori come in continuous families, but just in quadruples.
We shall also observe that considering planar 
 configurations in $\R^2\times \{(0,0)\}$ or in $\{(0,0)\}\times \R^2$ does not lead to anything interesting, since for such 
 configurations \eqref{BCnormalized} reduces to the (normalized) central configurations equation, 
 whereas considering any other plane spanned by $\{v_1,v_2\}$, where $v_1\in \R^2\times \{(0,0)\}$ 
 and $v_2 \in \{(0,0)\}\times \R^2$, do not produce different $\sbc$ by the $S^1\times S^1$-invariance of 
 \eqref{BCnormalized}. 
 
 
 Without further mentioning it, we will hereafter only consider $\sbc$ which are contained in $\{0\}\times \R^2\times \{0\}\cong \R^2$ (in other words, we will consider only solutions of \eqref{BCplane}), and refer to them simply as $\sbc$. Starting point for the results in \cite{AP:20}, and for the results of the present paper as well, is a careful study of the inertia indices of collinear $\sbc$, in shorthand notation $\csbc$. As one readily sees from \eqref{BCplane}, $\csbc$ must be contained in one of the two coordinate axes: we will 
henceforth call $s-\csbc$ those $\csbc$ which are contained in $\R\times \{0\}\subset \R^2$, and $1-\csbc$ 
those $\csbc$ which are contained in $\{0\}\times \R\subset \R^2$. 

A straightforward computation shows that, after rearranging properly the coordinates of the $s-\csbc$ $q$, we have the following block decomposition of the Hessian matrix 
$$H(q) = \left (\begin{matrix} -2 M^{-1}B(q) & 0 \\ 0 & M^{-1}B(q) \end{matrix}\right ) + \left (\begin{matrix} s U(q) I_n & 0 \\ 0 & U(q) I_n\end{matrix} \right ),$$
where $B(q)$ is the $(n\times n)$-matrix whose $(i,j)$-th entry is given by 
$$b_{ij}(q) = \frac{m_im_j}{r_{ij}^3}, \quad b_{ii} = - \sum_{j\neq i} \frac{m_im_j}{r_{ij}^3}.$$
Similarly, we have the following block decomposition of the Hessian matrix at any $1-\csbc$ $\widehat q$: 
$$H(\widehat q) = \left (\begin{matrix} -2 M^{-1}B(\widehat q) & 0 \\ 0 & M^{-1}B(\widehat q) \end{matrix}\right ) + \left (\begin{matrix}  U(\widehat q) I_n & 0 \\ 0 & s U(\widehat q) I_n\end{matrix} \right ).$$
Finally, we shall notice that $1-\csbc$ are actually normalized collinear central configurations, whereas $s-\csbc$ are obtained by scaling normalized collinear central configurations by a factor $1/\sqrt s$.
Putting these facts together, we proved in \cite[Section 2.2]{AP:20} the following result about the 
inertia indices of $\csbc$. In what follows we denote by $\iota^0(\widehat q),\ \iota^-(\widehat q),\ \iota^+(\widehat q)$ the nullity, Morse index, and Morse coindex respectively of a $1-\csbc$ $\widehat q$, with 
$$\eta_k(\widehat q) < ... < \eta_1(\widehat q)< \eta_0(\widehat q) := - U(\widehat q) <0$$
the distinct eigenvalues of the matrix $M^{-1}B(\widehat q)$, where $\widehat q$ is a fixed $1-\csbc$, 
and by $\alpha_k(\widehat q), ... , \alpha_1(\widehat q), \alpha_0 (\widehat q)= 1$ the corresponding multiplicities. Even if the eigenvalues of $M^{-1}B(\widehat q)$ and their multiplicities depend on the choice of $\widehat q$ in general, for the sake of readability we will hereafter drop the dependence on $\widehat q$.

\begin{prop}\label{thm:spectrum-csbc}
For any $s>1$, the inertia indices of any $s-\csbc$ $q$ are given by 
$$\iota^0 (\widehat q) = 0, \quad \iota^+(\widehat q) = n-2, \quad \iota^-(\widehat q) = n-1.$$
For any $1-\csbc$ $\widehat q$ we have:
\begin{enumerate}
\item If $\displaystyle -\frac{\eta_j}{U(\widehat q)} < s< - \frac{\eta_{j+1}}{U(\widehat q)}$, for some $j\in \{0,...,k-1\}$, then 
$$\iota^0 (\widehat q) = 0, \quad \iota^+(\widehat q) = n-2 + \sum_{i=0}^j \alpha_i, \quad \iota^-(\widehat q) = n-1 - \sum_{i=0}^j \alpha_i.$$
\item If $s = \displaystyle -\frac{\eta_j}{U(\widehat q)} $ for some $j\in \{1,...,k\}$, then 
$$\iota^0 (\widehat q) = \alpha_j, \quad \iota^+(\widehat q) = n-2 + \sum_{i=0}^{j-1} \alpha_i, \quad \iota^-(\widehat q) = n-1 - \sum_{i=0}^j \alpha_i.$$
In particular, $\widehat q$ is a degenerate critical point of $U$.
\item If $ s > - \displaystyle \frac{\eta_k}{U(\widehat q)}$, then 
$$\iota^0 (\widehat q) = 0, \quad \iota^+(\widehat q) = 2n-3, \quad \iota^-(\widehat q) = 0.$$
In particular, $\widehat q$ is a local minimum of $U$. \qed
\end{enumerate}
\end{prop}

From the proposition above we can easily deduce several facts: First, $\csbc$ are generically non-degenerate. Second, the inertia indices of $s-\csbc$ do not depend on $s$, whereas those of $1-\csbc$ strongly do. 
Even more, the Morse index of a $1-\csbc$ $\widehat q$ jumps at precisely characterized values of $s$ which only
depend on the spectrum of the matrix $M^{-1}B(\widehat q)$ and on the value of the Newtonian potential at $\widehat q$. This will enable us in the next section to show the existence of bifurcations of critical points of $\widehat U$ starting from $1-\csbc$. 

We shall also notice that in general we have no information about the eigenvalues of the matrix $M^{-1}B(\widehat q)$ and their multiplicities, besides the fact that $-U(\widehat q)$ is the 
largest non-zero eigenvalue and that it is simple. However, it is reasonable to believe that, for generic 
choice of the masses $m_1,...,m_n>0$, all eigenvalues of $M^{-1}B(\widehat q)$ are simple for any $1-\csbc$ $\widehat q$. 

Using Proposition~\ref{thm:spectrum-csbc} and the classical Morse inequalities, in \cite[Section 4]{AP:20} we gave the following lower bounds on the number of non collinear $\sbc$ assuming that all $\sbc$ are non-degenerate. 

\begin{thm}
\label{thm:lowerbound1}
Assuming that all $\sbc$ are non-degenerate, the following lower bounds hold:
\begin{enumerate}
\item If $\displaystyle s > \max \Big \{ - \frac{\eta_k(\widehat q)}{U(\widehat q)} \ \Big |\ \widehat q \textrm{ is  a }\ 1-\csbc\Big \}$, then there are at least 
$$3n! - 2(n-1)! - 2$$
non collinear $\sbc$. 
\item In all other cases, there are at least 
$$n! - 2(n-1)!$$
non-collinear $\sbc$. 
\end{enumerate}
As a corollary, for every $m_1,...,m_n>0$ and $s>1$ such that all $\sbc$ are non-degenerate, we have at least $n!-2(n-1)!$ relative equilibria of the form 
$$q(t) = \left (\begin{matrix} e^{i\sqrt{s} t} & 0 \\ 0 & e^{it}\end{matrix} \right )\cdot q, \quad q \textrm{ is a } \sbc,$$ 
for the $n$-body problem in $\R^4$ which are not induced by central configurations. \qed
\end{thm}

The lower bound given in Part 2 can be significantly improved by considering different cases and implementing the asymptotic estimates on the coefficients of the Poincar\'e polynomial of $\widehat{\mathbb S}$ proved in \cite[Section 3]{AP:20}. We refrain to do it here to keep the exposition as 
simple as possible. 

A celebrated result of Moeckel, known as the $45^\circ$-{\sc Theorem}, states that for $s=1$ (which we recall, corresponds to the central configurations case) the manifold of 
configurations which are collinear along some line is an attractor for the gradient flow of $U$ restricted to $\widehat {\mathbb S}$; in particular,
collinear central configurations are isolated, and an isolating set is given by the space of 
configurations for which the ``collinearity angle'' is less or equal to $45^\circ$. This can be seen as a global 
version of the fact that the Morse index of collinear central configurations is always as large as possible, namely $(d-1)(n-2)$ for the $n$-body problem in $\R^d$. In \cite[Section 5]{AP:20} we generalized the $45^\circ$-theorem to $\sbc$, proving (in the setting considered in the present paper) that the manifold of configurations which are collinear along the $x$-axis is an attractor for the gradient flow of $U$ restricted to $\widehat{\mathbb S}$. For the general $45^\circ$-theorem for $\sbc$ we refer to \cite[Theorem 5.6]{AP:20}.

\begin{thm}[$45^\circ$-theorem for $\sbc$]
\label{thm:45gradi}
The manifold 
$$\Big \{q=(q_1,...,q_n)\in \mathcal P \ |\ q_k \in \R\times \{0\}, \ \forall k =1,...,n \Big \}\subset \mathcal P$$
is an attractor for the gradient flow of $U$. More precisely, the {\sc collinearity function} 
$$\theta (q) := \max_{i\neq j}\  \angle (q_i-q_j, \partial_x )$$
is a Lyapounov function for the gradient vector field of $U$ on the set $\{q \in\mathcal P \ |\ 0<\theta(q) \leq 45^\circ\}.$\qed 
\end{thm}

We shall notice that a verbatim generalization of Moeckel's $45^\circ$-theorem to $\sbc$ is not possible since $1-\csbc$ are, for suitable values of $s>1$, local minima of $U$, and actually there is absolutely no reason why the manifold of configurations which are collinear along some line in $\R^2$ should be invariant under the gradient flow of $U$. In Section \ref{sec:4} we will strengthen these observations by showing that, for increasing value of $s>1$, we can find families of critical points of $U$ bifurcating from 
the set of $1-\csbc$. 

We finish this section recalling that the $45^\circ$-theorem for $\sbc$ can be used to improve the lower 
bounds given in Theorem \ref{thm:lowerbound1} on the number of relative equilibria in $\R^4$ assuming non-degeneracy; for more details we refer to \cite[Section 6]{AP:20}. We shall notice that here, unlike in Theorem \ref{thm:lowerbound1}, we are not able 
to exclude that such relative equilibria come from central configurations. Nevertheless, the result is 
still noteworthy since the lower bound that we obtain is larger than the largest known lower bound on the 
number of planar central configurations, see \cite{McCord96}.

\begin{thm}
For $s>1$ fixed, if all $\sbc$ are non-degenerate there are at least 
$$n! \left ( 1 + \frac 1n + \frac 32 \sum_{j=3}^n \frac 1j \right )$$
relative equilibria of the form 
$$q(t) = \left (\begin{matrix} e^{i\sqrt{s}} t & 0 \\ 0 & e^{it}\end{matrix}\right ) \cdot q, \quad q \textrm{ is a } \sbc,$$
for the $n$-body problem in $\R^4$. \qed
\label{thm:lowerbound2}
\end{thm}


\section{A brief recap on the spectral flow in finite dimension}\label{sec:sf}

The spectral flow is an integer-valued homotopy invariant of paths of selfadjoint Fredholm operators introduced by Atiyah, Patodi and Singer in the seventies in connection with the \textit{eta-invariant} and \textit{spectral asymmetry}.  In this section, we briefly recall the definition and the basic properties 
of the spectral flow in a finite dimensional setting. An elementary and self-contained introduction to the spectral flow for bounded selfadjoint Fredholm operators in infinite dimensional real Hilbert spaces can be found in \cite{FPR99}, whilst a quick recap and description of the spectral flow in the more general setting of paths of selfadjoint unbounded Fredholm operators having fixed domain appears in the beautiful paper \cite{RS95}. For further approaches we refer the interested reader to  \cite{HP17, HP19, HPY20} and references therein. 

Let $(H, \langle \cdot, \cdot \rangle)$ be an Euclidean  space and denote by $\mathcal L_{sym}(H) $ the vector space of all linear maps $T: H \to H$ that are self-adjoint with respect to $\langle \cdot, \cdot \rangle$. Roughly speaking, the {\sc spectral flow\/} $\spfl(L_t, t \in [a,b])$ of a continuous path $L: [a,b] \to \mathcal L_{sym}(H)$ is the number of negative eigenvalues of $L_a$ that become positive minus the number of positive eigenvalues of $L_a$ that become negative as the parameter $t$ runs from $a$ to $b$. In other words, the spectral flow measures the net change of eigenvalues crossing $0$ and can be interpreted as a sort of generalized signature. This informal description can be made rigorous in very many different ways. 

\begin{dfn}\label{def:sf}
Let $L:[a,b] \to 	\mathcal L_{sym}(H)$ be a continuous path of self-adjoint operators having invertible endpoints. We term {\sc spectral flow of $L$\/} on the interval $[a,b]$ the integer 
\[
\spfl(L_t, t \in [a,b]) \= \iota^-(L_a)-\iota^-(L_b)
\] 
where $\iota^-$ denotes the number of negative eigenvalues. \end{dfn}

A path of operators having invertible ends will be usually referred to as {\sc admissible.\/}
Under the non-degeneracy assumption on the endpoints (which is, for different reasons, 
always assumed throughout the paper) the (RHS) in the equation defining the spectral flow, can be 
equivalently written as 
\[
\dfrac12[\sgn(L_a)-\sgn(L_b)],
\]
thus pointing out why the spectral flow can be though of as a generalized signature. 
In this respect we observe that, assuming more regularity on the path $L$ (for instance, that $L$ is at least $\mathscr C^1$) it is possible to prove that the local contribution to the spectral flow is provided by the signature of a quadratic form (usually called {\sc crossing form\/}). More precisely, if $t_0 \in (a,b)$ is a {\sc crossing instant\/}, meaning that $\ker L(t_0)\neq \{0\}$, and if the restriction of the derivative of the path onto $\ker L(t_0)$   is non-degenerate as a quadratic form, then the spectral flow across the instant $t_0$ can be computed as follows 
\[
\spfl(L_t, t \in [t_0-\delta, t_0+\delta]) = \sgn \ L'(t_0)\vert_{\ker L(t_0)}.
\]

Here below we list some properties of the spectral flow that will be used in this paper. In what follows, 
every path of self-adjoint operators is assumed to be continuous. 

\begin{enumerate}
	\item {\sc Normalization.\/} Let $L:[a,b] \to \mathrm{GL}_{sym}(H)$ be a path of invertible operators. Then 
	\[
	\spfl(L_t, t \in [a,b])=0.
	\]
	\item {\sc Invariance under Cogredience.\/} If $L:[a,b] \to 	\mathcal L_{sym}(H)$ is admissible, then for any  $M: [a,b] \to \mathrm{GL}_{sym}(H)$ we have
	\[
	\spfl(L_t, t \in [a,b])=\spfl(M^*_tL_tM_t, t \in [a,b]). 
	\]
	\item {\sc Concatenation.\/} For $c\in[a,b]$, if $L:[a,b] \to \mathcal L_{sym}(H)$ is admissible on both $[a,c]$ and $[c,b]$, then 
	\[
	\spfl(L_t, t \in [a,b])= \spfl(L_t, t \in [a,c]) + \spfl(L_t, t \in[c,b])
	\]
	\item{\sc Homotopy invariance property.\/} If $H:[0,1]\times [a,b] \to \mathcal L_{sym}(H)$ is such that  the path $t \mapsto H(s, t)$ is admissible for each $s \in [0,1]$, then 
	\[
	\spfl(H(s,t), t \in [a,b]) = \spfl (H(0,t), t\in [a,b]),\quad \forall s \in [0,1].
	\]
\end{enumerate}

We finish this section observing that in Definition~\ref{def:sf} we required the path $L$ to have invertible endpoints. This assumption can be removed by properly choosing the contribution of the endpoints. In this more general setting, the spectral flow will be homotopy invariant with respect to end-points or, more precisely, end-points stratum homotopy invariant, meaning that the end-points are allowed to vary without changing the nullity. 


\section{Bifurcations of collinear $S$-balanced configurations}
\label{sec:4}

In this section, we prove an abstract bifurcation result from the trivial branch of a one parameter $\mathscr C^2$-family of functions on a finite-dimensional manifold and then apply it to the study of bifurcations $S$-balanced configurations. We start by introducing the natural geometric framework and some preliminary definitions. 

\begin{dfn}
A {\sc smooth family of finite-dimensional real smooth manifolds\/} $(X_\lambda)_{\lambda \in I}$ parameterized by the interval $I \=[a,b]\subset \R$ is a family of manifolds of the form $X_\lambda=p^{-1}(\lambda)$ where $p: X \to I$ is a smooth submersion of a manifold $X$ onto $I$.
\end{dfn}

For a smooth family $(X_\lambda)$ as above, $X_\lambda$ is a codimensional one submanifold of $X$ for every $\lambda \in I$. For each $x \in X_\lambda$, we have that $T_x X_\lambda = \ker D p_x$ and $$T^vX\=\Set{\ker Dp_x| x \in X}$$ is a vector subbundle of the tangent bundle $TX$. 

A smooth function $F:X \to \R$ defines a smooth family of functions $F_\lambda: X_\lambda \to \R$ by restriction to the fibers of $p$. We assume that there exists a smooth section $\sigma:I \to X$ of $p$ such that, for every $\lambda \in I$, $\sigma(\lambda)$ is a critical point of $F_\lambda$, and in what follows we refer to such a $\sigma$ as a {\sc trivial branch of critical points.\/}
\begin{dfn}\label{def:bfurcation-point}
	 We term $\lambda_* \in I$ a {\sc bifurcation instant from the trivial branch $\sigma(I)$\/} if there exists a sequence $\lambda_n \to \lambda_*$ and a sequence $(x_n)_{n \in \N} \subset X$ converging to $\sigma(\lambda_*)$ such that $p(x_n)=\lambda_n$ and each $x_n$ is a critical point for $F_{\lambda_n}$ not belonging to $\sigma(I)$.
\end{dfn}
\begin{note}
We hereafter denote the Hessian of $F_\lambda$ at the critical point $\sigma(\lambda)$ by $h_\lambda$.	
\end{note}

The family of Hessians $(h_\lambda)_{\lambda \in I}$ defines a smooth function $h$ on the
total space of the pull-back bundle $\mathcal H= \sigma^*(T^vX)$ of the vertical tangent
bundle. Hence, the function $h: \mathcal H \to \R$ defines a family of (generalized)
quadratic forms $h_\lambda$ defined on $T^v X_\lambda$.  Using the notation introduced in
Section \ref{sec:sf}, we say that the path $\lambda \mapsto h_\lambda$ is admissible if $h_a, h_b$ are non-degenerate. Denoting by $\spfl(h_\lambda, \lambda \in [a,b])$ the spectral flow of the path $h$, we have the following 
\begin{thm}\label{thm:main-abstract}
If $h=(h_\lambda)_{\lambda \in [a,b]}$ is admissible, then there exists at least one bifurcation instant $\lambda_*\in (a,b)$ of critical points of $F$ from the trivial branch.
Moreover, if $\ker h_\lambda \neq \{0\}$ only for finitely many $\lambda$, then there are at least 
\[
\dfrac{\big|\spfl(h_\lambda, \lambda \in [a,b])\big|}{m}
\] 
distinct bifurcation instants in $(a,b)$ where $m\= \max\ \{\dim \ker h_\lambda\}$. 
\end{thm}

\begin{rem}
The  {\sc singular set\/} $\Sigma(h)\=\Set{\lambda \in I | \ker h_\lambda \neq 0}$ is finite as soon as the data are analytic, which often occurs in the applications and is indeed the case in the situation considered in the present paper. \qed
\end{rem}

\begin{proof}
We split the proof into three steps: 

\vspace{2mm}

{\bf Step 1: (Reduction to a fixed Euclidean space)}
	By the vector bundle neighborhood theorem there exists  a trivial Euclidean bundle $\mathcal E=I \times H$  over $I=[a,b]$ and a fiber preserving smooth map $\psi: \mathcal E \to X$ such that $\psi(\lambda,0)=\sigma(\lambda)$ for every $\lambda \in [a,b]$, and $\psi$ is a diffeomorphism of $\mathcal E$ into an open neighborhood $\mathcal O$ of $\sigma(I)$ in $X$. 

Let $\widetilde F: I \times H \to \R$ be the map defined by composition $\widetilde F \= F \circ \psi$. So, the map $\widetilde F$ defines a smooth one parameter family of functions on $H$. Since $\psi$ is a fiber preserving diffeomorphism, $u \in H$ is a critical point of $\widetilde F_\lambda$ if and only if $x = \psi_\lambda(u)$  is a critical point of $F_\lambda$. In particular, $0$ is a critical point of $\widetilde F_\lambda$ for each $\lambda \in I$. 

The Hessian $\widetilde h_\lambda$ of $\widetilde F_\lambda$ at $0$ is given by $\widetilde h_\lambda(\zeta)=h_\lambda(D_0\psi_\lambda)([\zeta])$. By the cogredience and normalization properties of the spectral flow, see Section 3, we get that 
\[
\spfl(\widetilde h_\lambda, \lambda \in I)=\spfl(h_\lambda, \lambda \in I)= \spfl(L_\lambda, \lambda \in I)= \iota^-(L_a)- \iota^-(L_b)
\] 
where $L: I \to  L^s(H)$ is a smooth path of self-adjoint operators representing the quadratic form $\widetilde h$ with respect to the scalar product of $H$, namely $\widetilde h_\lambda(u)=\langle L_\lambda u, u\rangle$ for every $u \in H$.

\vspace{2mm}

{\bf Step 2: (Non-vanishing spectral flow implies bifurcation)} 
We assume by contradiction that 
$$\spfl(L_\lambda, \lambda \in I)= \iota^-(L_a)- \iota^-(L_b)\neq 0 $$ 
and that there are no non-trivial critical points of the path $\lambda \mapsto \widetilde F_\lambda$ close to the trivial branch. Then, there exists $\delta >0$ such that $0\in H$ is the only critical point of $\widetilde F_\lambda$ on $B_\delta$ for every $\lambda \in I$, where with $B_\delta \subset H$ we 
denote the open ball with radius $\delta$ around the origin. Without loss of generality we can suppose $\widetilde F_\lambda (0)=0$ for every $\lambda\in I$.
For any $\lambda \in I$ and any non-negative integer $k$,  let $C_k(\widetilde F_\lambda, 0)$ be the $k$-th local homology group associated to the isolated critical point $0$ of the functional $\widetilde F_\lambda$ 
\[
C_k(\widetilde F_\lambda, 0)\= H_k(\widetilde F_\lambda ^{0} \cap B_\delta,\widetilde F_\lambda ^{0} \cap B_\delta \setminus\{0\})
\]
where, as usually, $H_k(\cdot, \cdot)$ denote the $k$-th relative singular homology group with integer coefficients, and
$$\widetilde F_\lambda ^{0}\=\Set{x \in H | \widetilde F_\lambda(x) \le 0}$$
denotes the sublevel set of $\widetilde F_\lambda$. Since by assumption $0$ is an isolated critical point 
of $F_\lambda$ for every $\lambda \in I$, for each $k\in \N_0$ the rank of the $k$-th local homology group is independent of $\lambda$. Moreover, the admissibility of the path implies in virtue of the Morse Lemma that
\[
C_k(\widetilde F_a, 0)= \left \{
\begin{array}{r}
	\Z  \quad \textrm{ if } k= \iota^-(L_a),\\
	0  \ \ \, \quad \quad  \textrm{ otherwise,}
\end{array}\right .  \qquad \qquad 
C_k(\widetilde F_b, 0)=\left \{
\begin{array}{r}
	\Z  \quad \textrm{ if } k= \iota^-(L_b).\\
	0  \ \ \quad \quad \textrm{ otherwise.}
\end{array}\right .
\]
However, since the spectral flow on the interval $I$ does not vanish, we have that $\iota^-(L_a) \neq \iota^-(L_b)$, in contradiction with the fact that the local homology groups do not depend on $\lambda$.

\vspace{2mm}

{\bf Step 3: (Estimate from below on the number of bifurcation instants)} In order to prove this last assertion we proceed as follows. By assumption, there are only finite many instants $a < \lambda_1 < \ldots < \lambda_k < b$ for which the kernel of $L_\lambda$ is non-zero. Therefore, we can find $\delta>0$ such that for every $j =1, \ldots, k$, the intervals $I_j\= [\lambda_j-\delta, \lambda_j+\delta]$ are pairwise disjoint and $L_\lambda$ is non-degenerate at the endpoints of all such intervals.  By the additivity property of the spectral flow, we get that 
\[
\spfl(L_\lambda, \lambda \in I)=\sum_{j=1}^k \spfl(L_\lambda, \lambda \in I_j)=\sum_{j=1}^k \big[\iota^-(L_{\lambda_j-\delta})- \iota^-(L_{\lambda_j+\delta})\big]
\]  
Since $\dim \ker L_\lambda \le m$ it follows that 
\[
\big\vert\iota^-(L_{\lambda_j-\delta})- \iota^-(L_{\lambda_j+\delta})\big\vert\le m.
\]
Summing up, we immediately get that 
\[
\big\vert\spfl(L_\lambda, \lambda \in I)\big\vert \le \sum_{j=1}^k \big\vert \spfl(L_\lambda, \lambda \in I_j)\big\vert \le dm
\]
where $d$ denotes the number of non-vanishing terms in the sum. Since by Step 2 any interval $I_j$ having non-vanishing spectral flow contributes to the total number of bifurcation instants at least by 1, we immediately get that there must exist at least $d= |\spfl(L_\lambda, \lambda \in I)|/m$ bifurcation points on the interval $I$. 
\end{proof}

\begin{rem}
A similar result holds in infinite dimensional separable Hilbert spaces for a path of (bounded) self-adjoint Fredholm operators even in the strongly indefinite case, meaning that the Fredholm operators are compact perturbations of  an invertible operator having an infinite dimensional positive and negative eigenspace. In this case in fact both the Morse index and coindex are infinite. For further details we refer the interested reader to \cite{MPP07, PPT04, PW14, PW17} and references therein for this more general functional analytic setting with application e.g. to semi-Riemannian geodesics and conjugate points. \qed
\end{rem}

Theorem~\ref{thm:main-abstract} will be the key ingredient for our bifurcation result for $\sbc$, where the bifurcation parameter is precisely the parameter $s$ appearing in the matrix $S$. Thus, for every $s>1$
we consider the function $U : \mathcal P_s \to \R$. By the discussion provided in Section~\ref{sec:preliminaries},  $\sbc$ and in particular $\csbc$ are critical points of $U$. The latter turn out to 
be generically non-degenerate, see Proposition~\ref{thm:spectrum-csbc}, meaning that the Hessian of $U$ at any $\csbc$ is non-degenerate as a quadratic form for all but finitely many values of $s$. Before stating the 
bifurcation theorem for $\sbc$, we shall finally observe that $1-\csbc$, which we recall are nothing else but normalized collinear central configurations along the $y$-axis, are solutions of \eqref{BCplane}
independently of $s >1$. In particular, any constant one-parameter family $(\widehat q_s = \widehat q)_{s>1}$, with $\widehat q$ a fixed $1-\csbc$, provides a trivial branch of critical points. 

\begin{thm}\label{main:sbalanced-bifurcation}
For $s_1, s_2 >1$, let $J\=[s_1, s_2]$ and let $(\widehat q_s)_{s \in J}$ be the trivial one-parameter family of $\sbc$ given by a fixed $1-\csbc$ $\widehat q$, that is $\widehat q_s = \widehat q$ for every $s\in J$. Setting 
$$\alpha_*(\widehat q)\= \max\Set{\alpha_j(\widehat q) | j =1, \ldots, k },$$ 
for  $s_1$ sufficiently close to $1$ and $s_2$ large enough there exist at least 
\[
\left \lfloor \frac{n-2}{\alpha_*(\widehat q)}\right \rfloor 
\]
bifurcation instants from $\widehat q$. As a corollary, for  $s_1$ sufficiently close to $1$ and $s_2$ large enough there are at least 
$$n! \cdot \left \lfloor \frac{n-2}{\alpha_*}\right \rfloor, \quad \alpha_* := \max \{\alpha_*(\widehat q) \ |\ \widehat q \textrm{ is a } 1-\csbc\},$$ 
bifurcation instants from the trivial families of $1-\csbc$. 
\end{thm}

\begin{rmk}
 Even if we have no information about the $\alpha_j$'s in general, it is reasonable to believe that they are
all equal one for generic choice of the masses $m_1,...,m_n$. In this case, 
Theorem~\ref{main:sbalanced-bifurcation} would provide the existence of at least $n! (n-2)$ bifurcation instants, namely at least $n-2$ bifurcation instants for each choice of the ordering of the masses along the $y$-axis.  \qed
\end{rmk}

\begin{proof}
We denote by  $\pi: \mathcal P \to J$ the (trivial) ellipsoid bundle (topologically, a sphere bundle) over $J$ whose fiber $\pi^{-1}(\{s\})= \mathcal P_s $ is the collision free configuration sphere. Notice that $\mathcal P_s$  depends on the parameter $s$; in particular, for $s=1$ we obtain a round sphere in the mass metric, which for increasing $s$ is deformed into an ellipsoid having its major axes in the directions corresponding to the eigenvalue $1$ of the matrix $S$. 

The Newtonian potential $U$ defines a smooth bundle map $\mathcal U$ whose restriction $\mathcal U_s$ to each fiber $\pi^{-1}(\{s\})$ is precisely $U$, and for each $s \in J$ 
the configuration vector $\widehat q_s=\widehat q$ is a $1-\csbc$, hence a critical point of $\mathcal U_s$. By Proposition~\ref{thm:spectrum-csbc}, we infer that the path $s\mapsto \widehat q_s$ is admissible 
(meaning that the associated path $s\mapsto h_s$ of quadratic forms pointwise defined by the Hessian of $U$ at $\widehat q_s$ is admissible) as 
soon as 
$$1<s_1 < -\frac{\eta_1(\widehat q)}{U(\widehat q)}, \qquad s_2 > -\frac{\eta_k(\widehat q)}{U(\widehat q)}.$$ 
Moreover, by setting $j=0$ in Proposition~\ref{thm:spectrum-csbc}, Part 1, we get that the $\iota^-(\widehat q_{s_1})= n-2$, whereas using Item iii) of Proposition~\ref{thm:spectrum-csbc}, Part 2, we get that 
$\iota^-(\widehat q_{s_2})= 0$. 
In particular, the spectral flow of the path $s\mapsto h_s$ is easily computed to be
\[
\spfl(h_s, s \in [s_1, s_2])= \iota^-(h_{s_1})- \iota^-(h_{s_2})= n-2
\]
where $ \iota^-(h_{s_1})\= \iota^-(\xi_{s_1})$ and $\iota^-(h_{s_2})\=\iota^-(\xi_{s_2})$.  Theorem~\ref{thm:main-abstract} now implies the claim, observing that the kernel of $h_s$ is non-trivial only in correspondence of the eigenvalues $\eta_j(\widehat q)$ whose multiplicity is $\alpha_j(\widehat q)$.
\end{proof}


\section{Some explicit examples for $n=3$}
\label{sec:5}
In this section we use numerical computations to study the non-trivial branches bifurcating
from a trivial branch of $1-\csbc$. First, we introduce the continuation method used to compute curves of
solutions of a system of nonlinear equations depending on a real parameter, then we
describe how to use it to follow the bifurcations branching from a trivial branch of $1-\csbc$. Finally, we
provide some examples for $n=3$. 

\subsection{The continuation method}

Let $F:\R^{n} \times \R \to \R^n$ be a differentiable function, and denote with $q \in
\R^n$ the spatial component and with $s \in \R$ the parameter. Take $(q_i, s_i)^T \in \R^n
\times \R$ such that $F(q_i, s_i) = 0$: the purpose of a continuation method is to find a
zero of $F$ for a different value of the parameter $s$, starting from the known solution at
$s_i$. 

To this end, the pair $(q,s)^T$ is displaced by solving the equation
\begin{equation}
   \begin{cases}
      F(q, s) = 0, \\[2ex]
   \left\lvert 
   \begin{pmatrix}
      q \\ s
   \end{pmatrix}
   -
   \begin{pmatrix}
      q_i \\ s_i
   \end{pmatrix}
   \right\rvert^2 - \delta^2 = 0,
   \label{eq:pseudoArcLength}
   \end{cases}
\end{equation}
where $\delta> 0$ is small. A solution of Equation \eqref{eq:pseudoArcLength} can be computed using the Newton method
(see e.g. \cite{SB02}), thus solving at
each step a system of equations given by the matrix
\begin{equation}
   \begin{bmatrix}
      \displaystyle \frac{\partial F}{\partial q} & \displaystyle \frac{\partial
      F}{\partial s} \\[2ex]
      2(q-q_i) & 2(s-s_i)
   \end{bmatrix}.
   \label{eq:systNewton}
\end{equation}
A first guess $(\bar{q}, \bar{s})$ for the Newton method can be constructed starting from the known solution
$(q_i, s_i)^T$, and taking a tangent displacement along the curve of solutions. 
The tangent line can be approximated by using two consecutive solutions $(q_i, s_i)^T$ and
$(q_{i-1}, s_{i-1})^T$, hence the first guess can be taken as
\begin{equation}
   \begin{pmatrix}
      \bar{q} \\
      \bar{s}
   \end{pmatrix}
   =
   \begin{pmatrix}
      q_i \\ s_i
   \end{pmatrix}
   + \gamma
   \begin{pmatrix}
      q_i - q_{i-1} \\
      s_i - s_{i-1}
   \end{pmatrix},
   \qquad
   \gamma = \frac{\delta}{\lvert (q_i, s_i)^T - (q_{i-1}, s_{i-1})^T \rvert}.
   \label{eq:approxTangent}
\end{equation}
If not known, an additional solution $(q_{i-1}, s_{i-1})^T$ can be computed from the known
solution $(q_i, s_i)^T$ by simply displacing $s_{i}$ as $s_{i-1} = s_i + \Delta s$, and
then solving the equation $G(q) := F(q , s_{i-1}) = 0$ with the Newton method,
constructing a first guess using a displacement of $q_i$.

In the case of $S$-balanced configurations, the function $F$ reads
\[
   F(q,s) = M^{-1} \nabla U(q) + U(q) \widehat{S}(s) q,
\]
where $\widehat{S}(s)$ is the block-diagonal matrix
$$\widehat{S}(s) = \left ( \left ( \begin{matrix} s & 0 \\ 0 & 1 \end{matrix}\right ), ... , \left ( \begin{matrix} s & 0 \\ 0 & 1 \end{matrix}\right )\right ).$$

To numerically study the bifurcations branching from a trivial branch of $1-\csbc$ in the case of three masses, we use the following scheme. Let $m_1,m_2,m_3$ be given positive masses, and let an 
ordering of the three masses along the $y$-axis be fixed. Then,
\begin{enumerate}
   \item we compute the corresponding $1-\csbc$ (which we recall is nothing else but a normalized collinear central
      configuration), namely using the Newton method to find the unique zero of the Euler polynomial (see e.g. \cite{Moe94}) 
      corresponding to the chosen ordering of the masses. 
   \item We compute the eigenvalues of the matrix $M^{-1}B(\widehat{q})$, and the value of the
      potential $U(\widehat{q})$.
   \item We compute the (unique in virtue of Item ii) of Proposition \ref{thm:spectrum-csbc}) value $\tilde{s}$ for which $\widehat{q}$ is a degenerate critical point
      of $U$.
   \item We displace $\tilde s$ by $s = \tilde{s} + \Delta s$,
      and search for a zero of $F( \, \cdot \,, s)$ using the Newton method. A first guess is constructed by displacing
      $\widehat{q}$ along the direction of the kernel of $\partial F/\partial q$. 
   \item Using the previous two solutions, we start the continuation method with respect to the parameter $s$. 
\end{enumerate}
The same scheme can be used also for $n\geq 4$, modifying the first
step for the computation of the corresponding $1-\csbc$. Notice indeed that also in this case
there exists a unique $1-\csbc$ for each ordering of the bodies \cite{Mou10}. We plan to numerically investigate this case in future work.

\subsection{Results of the computations}
\begin{figure}[!ht]
   \centering
   \includegraphics[width=0.45\textwidth]{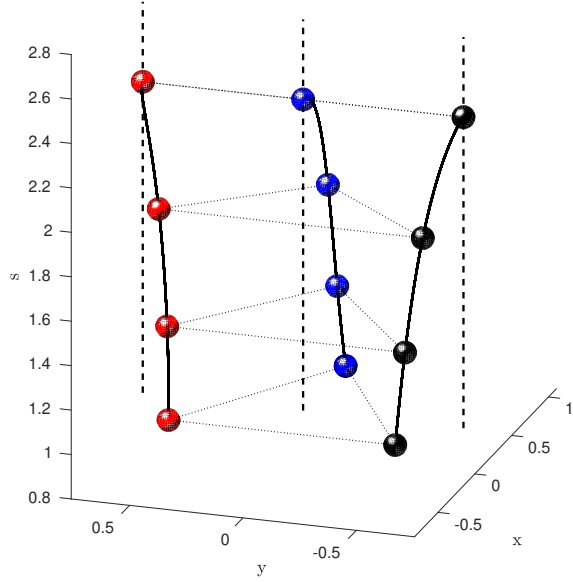}
   \caption{The branches of $\sbc$ in the case of three equal masses in the space
   $(x,y,s)$. The three dashed vertical lines represent the trivial branch of $1-\csbc$, 
   while the black and the red thick curves represent the two
   non-trivial branches originating from the bifurcation point.}
   \label{fig:equalMasses}
\end{figure}

We produced a first example using three unitary equal masses. The corresponding 
trivial branch of $1-\csbc$ has a bifurcation at $\tilde{s}=2.4$, from which two
non-trivial branches originate. On these branches the parameter $s$ decreases, while the three
masses move on isosceles configurations until arriving at an equilateral
configuration for $s=1$. The two branches are symmetrical with respect to the $y$-axis,
with the difference that the pair $\{ q_2-q_1, q_3-q_1 \}$ is positively oriented on one
branch and negatively oriented on the other one. Moreover, all solutions found are local
minima of the potential $U$. Figure~\ref{fig:equalMasses} shows one of the branches in the
space $(x,y,s)$; the other branch is not displayed for visibility reasons. Observe that, due to the additional symmetries of the problem, changing the ordering of the masses 
does not yield qualitatively different behaviors of the bifurcations branches. Also, a branch (actually four branches due to the $\Z_2\times \Z_2$-symmetry) of saddle points originating at $s=1$ from Lagrange's equilateral triangle is present but not shown in the figure. 

\begin{figure}[!ht]
   \centering
   \includegraphics[width=0.45\textwidth]{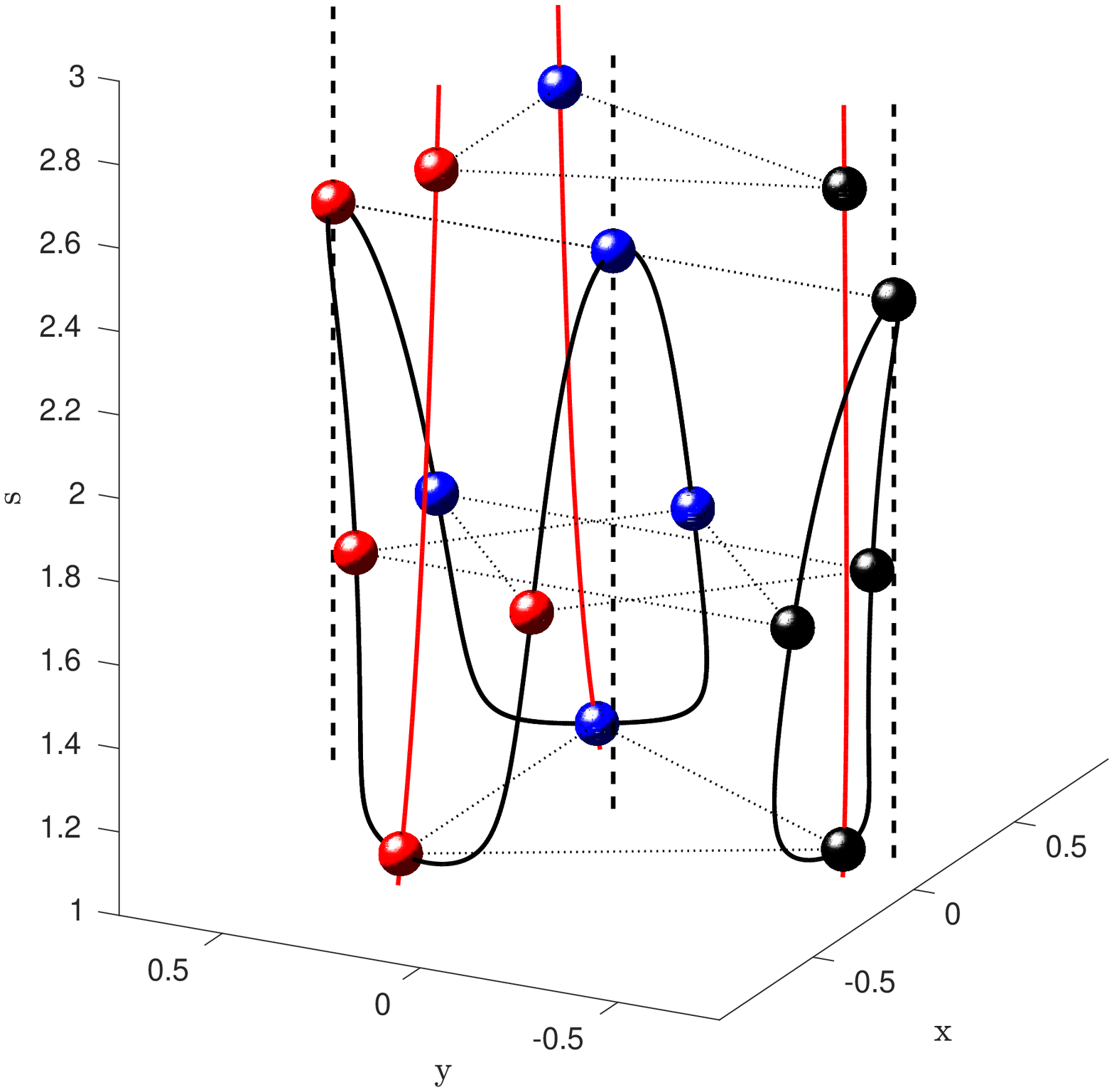}
   \includegraphics[width=0.45\textwidth]{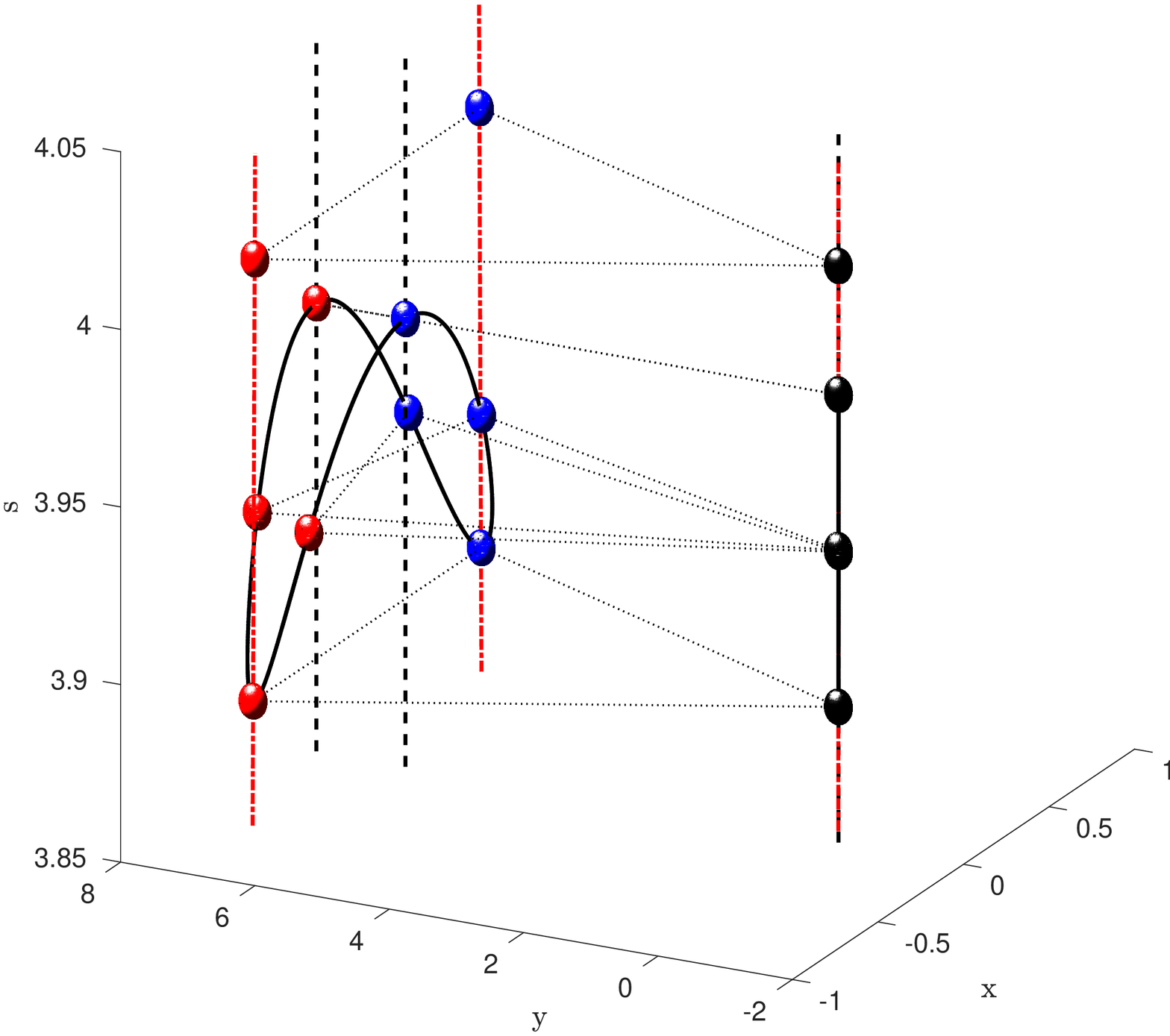}
   \caption{The branches of $\sbc$ in the case $m_1=1, \, m_2=m_3=\mu$, $0<\mu<1$,
   in the space $(x,y,s)$. The left panel refers to the case $\mu =0.99$, while the right
   one to $\mu = 0.01$. The black particle indicates the mass $m_1$, while the 
   blue and red particles refer to $m_2$ and $m_3$ respectively. As in
   Figure~\ref{fig:equalMasses}, the dashed vertical lines represent the trivial branch of $1-\csbc$,
   whereas the thick black curve represents a branch originating from the bifurcation
   point. Another branch, in which the masses follow the same curve in the opposite
   direction, is also present but not shown. The secondary branches 
   originating in correspondence of the turning point are drawn in red.}
   \label{fig:big1small2}
\end{figure}

We produced a second example using one unitary mass and two smaller equal masses,
say $m_1=1, \, m_2=m_3=\mu$, $0<\mu<1$, for the ordering of the masses given by taking the bigger mass outside the segment joining the two smaller ones. 
Also here we have a unique bifurcation point
$\tilde{s}$ from which two non-trivial (again, symmetrical with respect to the $y$-axis) 
branches originate. Along each non-trivial branch the parameter $s$ decreases at first until reaching a turning point $s_{\mathrm{turn}}$, where the three masses are 
placed on the vertices of a symmetric with respect to the $y$-axis isosceles triangle. After the turning point, the parameter $s$ increases up to the bifurcation value $\tilde s$, and the corresponding 
branch of $\sbc$ reaches for $s\uparrow \tilde s$ the $1-\csbc$ corresponding to the ordering of the masses in which the two smaller masses are swapped. 
In the $(x,y)$ plane, the masses $m_2$ and $m_3$ appear to rotate around a common point, while the mass $m_1$ appears to oscillate. 
The oscillation becomes larger and larger as the value $\mu$ approaches $1$, and the isosceles configuration at the turning point tends to an
equilateral one, which is reached in the limit $\mu\uparrow 1$. On the other hand, the oscillations are very small for $\mu \ll 1$,
making the position of $m_1$ almost constant. Also here, all configurations along the branches are local minima of the potential $U$. In correspondence of the turning point $s_{\mathrm{turn}}$
we also have two additional branches originating from the non-trivial branches, namely a first one along which the parameter $s$ decreases to 1 reaching in the limit $s\downarrow 1$ the 
equilateral configuration, and a second one along which the parameter $s$ increases to infinity reaching in the limit $s\to +\infty$ a limit configuration. The critical points 
along the first secondary branch are local minima of $U$, whereas along the latter one we have saddle points. 
The non-trivial branches in the space $(x,y,s)$ corresponding to the cases $\mu=0.99$ and $\mu=0.01$ respectively are depicted in 
Figure~\ref{fig:big1small2}. In the limit $\mu\uparrow 1$ the non-trivial branches ``tend'' to the non-trivial branches depicted in Figure~\ref{fig:equalMasses}.

In case the unitary mass is placed between the two smaller ones, the behavior is similar to the three equal masses case: the parameter
$s$ decreases down to $1$, where the masses reach an equilateral configuration.

\begin{figure}[!ht]
   \centering
   \includegraphics[width=0.45\textwidth]{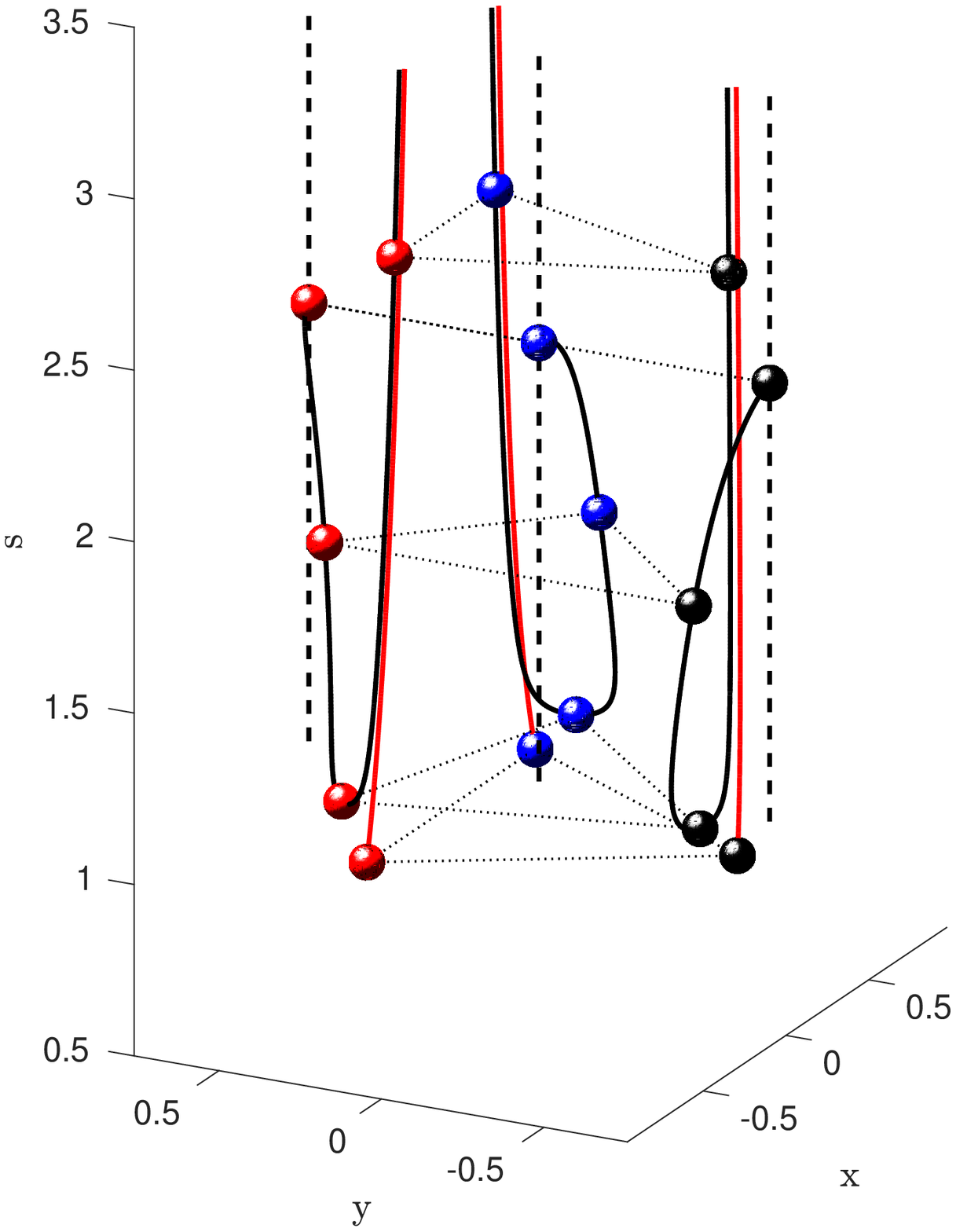}
   \includegraphics[width=0.45\textwidth]{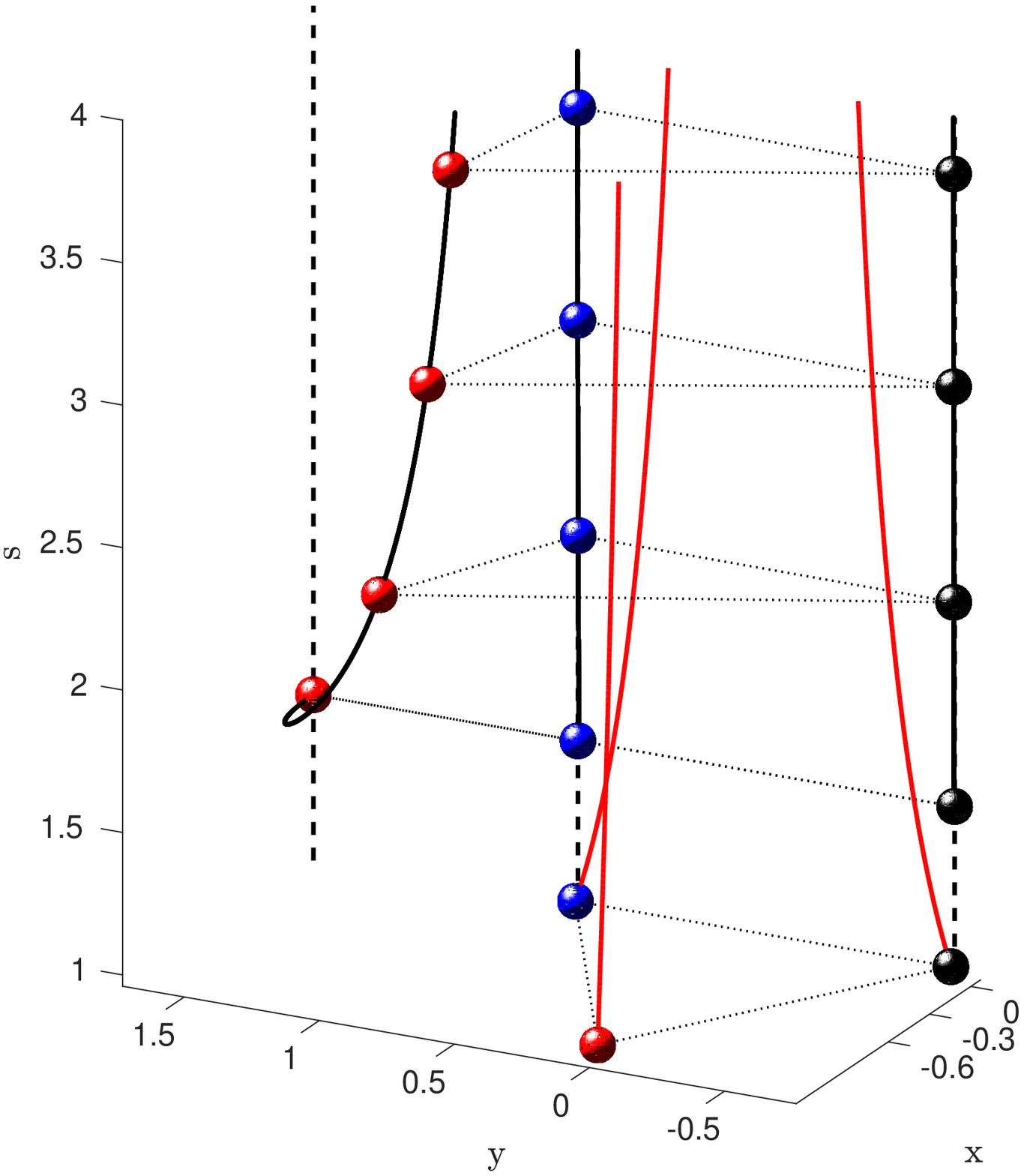}
   \caption{The branches of $\sbc$ in the case $m_1=m_2=1, \, m_3=\mu$, $0<\mu<1$,
   in the space $(x,y,s)$. The left panel refers to the case $\mu =0.99$, while the right
   panel to the case $\mu = 0.01$. The black particle indicates the mass $m_1$, while the 
   blue and the red particles refer to $m_2$ and $m_3$, respectively. As in
   Figures~\ref{fig:equalMasses} and \ref{fig:big1small2}, the dashed vertical line represents
   the trivial branch, while the thick black curve represent one branch originating from the bifurcation
   point. Another branch, symmetrical with respect to the $y$-axis, is also present but
   not shown. The branch of saddle points originating from the equilateral triangle is drawn in red.
   The black and red branches are asymptotic for $s\to +\infty$ to different asymptotic configurations.}
   \label{fig:big2small1}
\end{figure}

A third example has been produced using two unitary masses and a smaller one,
say $m_1=m_2=1, \, m_3=\mu, \, 0<\mu<1$, for the ordering of the masses in which the smaller mass lies outside the segment joining the two bigger ones. 
Two non-trivial branches, symmetrical with respect to the $y$-axis, originate from the (unique) bifurcation point $\tilde{s}$: on these
branches the parameter $s$ initially decreases, and the masses are placed on
scalene configurations corresponding to local minima of $U$.  A turning point $s_{\mathrm{turn}}$ is reached during the continuation, and the differential of $F$ is
singular at $s_{\mathrm{turn}}$. For $\mu\sim 1$, the turning point $s_{\mathrm{turn}}$ is close to $s=1$, and again the non-trivial branches approach the non-trivial branches depicted in Figure~\ref{fig:equalMasses} 
in the limit $\mu\uparrow 1$. On the other hand, the turning point becomes closer and closer to the bifurcation value $\tilde{s}$ for $\mu \to 0$.
After the turning point, the parameter $s$ increases to $+\infty$, with the masses
still placed on scalene configurations which in this case, however, correspond to saddle points of
$U$. Therefore, we have a jump on the Morse index in correspondence of the turning point $s_{\mathrm{turn}}$.
We also looked for secondary branches at the turning point $s_{\mathrm{turn}}$. To construct an initial guess for the Newton method, we displaced the configuration
at $s_{\mathrm{turn}}$ in two different manners: 
1) along the direction of the kernel of  $\partial F/\partial q$, and
2) displacing $q_1, q_2, q_3$ along random directions $v_1, v_2, v_3 \in S^1$,
respectively. In both cases, we did not find any other branch. Observe that this does not contradict Theorem~\ref{thm:main-abstract}, since the considered branch is not trivial. 
We plan to examine the occurrence of bifurcation points on non-trivial branches in future work.
Following the non-trivial branch after the turning point, we see that the configurations are asymptotic as $s$ grows to infinity to a limit configuration. As the left panel in Figure~\ref{fig:big2small1} might 
suggest, the asymptotic configuration seems to be the limit also of a branch of saddle points originating at $s=1$ from the equilateral triangle. This is in fact not the case, as it can be easily deduced 
from the right panel. 
The branches in the space $(x,y,s)$ corresponding to $\mu=0.99$ and $\mu=0.01$ respectively are depicted in 
Figure~\ref{fig:big2small1}.

In the case that the mass $m_3$ lies between $m_1$ and $m_2$, the result is similar to the three equal masses case, with the bifurcation parameter $s$ decreasing to $1$, where an
equilateral configuration is reached.

Further examples have been produced using three different masses by means of an additional parameter, but no qualitatively different behavior has been observed. 
More precisely, for fixed parameters $\mu,\lambda <1$, we considered masses $(m_1,m_2,m_3) = (1, \lambda, \mu)$. As it turns out, the qualitative behavior 
of the corresponding non-trivial branches only depends on $\lambda$: For $\lambda < \min\{1,\mu\}$ we are in situation analogous to the one depicted in Figure~\ref{fig:equalMasses}, with two 
branches originating from Euler's configuration on which $s$ decreases to $1$ reaching in the limit $s\downarrow 1$ Lagrange's equilateral triangle along minima of $U$. 
For $\lambda =\min\{1,\mu\}$, we have the same qualitative behavior as in Figure~\ref{fig:big1small2}: along each non-trivial branch the parameter $s$ decreases at first until reaching 
a turning point. After the turning point, the parameter $s$ increases up to the bifurcation value $\tilde s$, and the corresponding 
branch of $\sbc$ reaches for $s\uparrow \tilde s$ the $1-\csbc$ corresponding to the ordering of the masses in which the two smaller masses are swapped. In correspondence 
of the turning point two additional branches originate, namely one along which $s$ decreases reaching for $s=1$ Lagrange's equilateral triangle on minima of $U$, and one along which 
$s$ increases to infinity on saddle points of $U$. Finally, in the case $\lambda > \min \{\mu,1\}$ we have the same qualitative behavior as in
Figure~\ref{fig:big2small1}: The parameter $s$ 
decreases until a turning point is reached, after which
$s$ grows indefinitely and the configuration tends to an asymptotic configuration. 

Videos showing how the position of the masses in the $(x,y)$-plane changes along a non-trivial branch, for the cases discussed in this section, can be found at \cite{FenWeb}.
 
 
 \subsection{Final comments}
 
 We are now in position to make some final comments on the numerical implementations discussed in the previous subsection and on their consequences on the dynamics 
 of the $n$-body problem in $\R^4$. Before doing that, we shall recall that the bifurcations of $\sbc$ we found can be seen as planar $\sbc$ in $\R^4$ contained 
 in the plane $\{0\}\times \R^2\times \{0\}\subset \R^4$, and as such define relative equilibria for the $n$-body problem in $\R^4$ by 
 $$\mathbf{q}(t) = \left (\begin{matrix} e^{i \sqrt{s} t} & 0 \\ 0 & e^{it}\end{matrix}\right ) \cdot \left (\begin{matrix} 0 \\ x \\ y \\ 0\end{matrix}\right ),$$
 where $(0,x,y,0)^T=q$ is any planar $\sbc$ contained in $\{0\}\times \R^2\times \{0\}$.
 
 In Figure 1 we observe for the case of three equal masses a connection between Lagrange's equilateral triangle $q_{\Lag}$ and Euler's collinear configuration $q_\eul$ through a branch of $\sbc$
 which are local minima of $U$. This yields a continuum of relative equilibria $(\mathbf{q}_s)_{s \in [1,\tilde s]}$ for the $3$-body problem in $\R^4$
 $$\mathbf{q}_{s}(t) = \left (\begin{matrix} e^{i \sqrt{s} t} & 0 \\ 0 & e^{it}\end{matrix}\right ) \cdot \left (\begin{matrix} 0 \\ x_s \\ y_s \\ 0\end{matrix}\right ), \quad  \left (\begin{matrix} 0 \\ x_1 \\ y_1 \\ 0\end{matrix}\right ) = q_\Lag, \  \left (\begin{matrix} 0 \\ x_{\tilde s} \\ y_{\tilde s} \\ 0\end{matrix}\right ) = q_\eul,$$
 which is precisely the one found abstractly in \cite[Corollary 2.15]{AP:20}.
 An illustration of $\mathbf{q}_1$ can be found in \cite[Figure 9]{Moe:S}. 
 For $s>\tilde s$ such a continuum can be extended simply by taking $q_s =q_\eul$ for all $s>\tilde s$. Finally, for $s\in [1,\tilde s) \cap \Q$ or $s>\tilde s$, the relative 
 equilibrium $\mathbf{q}_s$ is periodic, and it is quasi-periodic otherwise. 
 
A similar connection can be observed also for arbitrary masses, namely between Lagrange's equilateral triangle and
Euler's collinear configuration in which the smaller mass lies between the other two masses. Such connections cannot be found 
within the class of central configurations, as Euler's configurations are isolated in virtue of Moeckel's $45^\circ$-theorem. 

As the implementations suggest, changes in the qualitative behavior of the non-trivial branches can be expected only while deforming the masses parameters by passing through a configuration
of the masses that forces additional symmetries of the problem, namely a configurations of the masses in which at least two of the masses are equal. We plan to investigate this aspect further for larger values of $n$. 

Another interesting group of open questions concerns the stability of relative equilibria generated by a $\sbc$. As already observed, to a given planar $\cc$ it is possible to associate a relative equilibrium solution in which all bodies rigidly rotate around their center of mass. In rotating coordinates, this relative equilibrium reduces to an equilibrium (actually, to the $\cc$ originating the orbit through a circle action).  So, it is quite natural to ask whether there is or not a relation between the linear stability of the relative equilibrium (as a periodic orbit) and the Morse index of the associated $\cc$ (in the rotating frame). Several results in this direction are nowadays available in the literature: For instance, in \cite{HS09} the authors provide a sufficient condition for the linear  instability of a relative equilibrium originated by a non-degenerate $\cc$ in the plane in terms of the Morse index, which in \cite{BJP14} is generalized to a broader class of singular operators as well as to the case of relative equilibria generated by a possibly degenerate critical point. 
In this latter case, besides the Morse index, a key role is played by the Jordan normal form associated to the Floquet multiplier 1. 

On the other hand, nothing seems to be known about the linear (and spectral) stability for higher-dimensional relative equilibria, see e.g. \cite{Moe14}. In higher-dimension, the situation is indeed much more involved:
For relative equilibria which are periodic in time (that is, for $s\in \Q$) we expect, besides Morse index and Jordan normal form, also the parameter $s$ to play a key role in the characterization of the stability properties. 
For relative equilibria which are quasi-periodic in time (that is, for $s\in \R\setminus \Q)$ a similar characterization relating the KAM stability of the relative equilibrium to the inertia indices of the associated $\sbc$
should be possible as well. All these questions will be addressed in future work.

\vspace{1cm}
\noindent
\textsc{Dr. Luca Asselle}\\
Justus Liebig Universit\"at Gie\ss en\\
Arndtrstrasse 2\\
35392, Gie\ss en\\
Germany\\
E-mail: $\mathrm{luca.asselle@math.uni}$-$\mathrm{giessen.de}$

\vspace{5mm}
\noindent
\textsc{Dr. Marco Fenucci}\\
Department of Astronomy, Faculty of Mathematics\\
University of Belgrade\\
Studentski trg 16 \\
11000 Belgrade \\
Serbia\\
E-mail: $\mathrm{marco\_fenucci@matf.bg.ac.rs}$

\vspace{5mm}
\noindent
\textsc{Prof. Alessandro Portaluri}\\
Università degli Studi di Torino\\
Largo Paolo Braccini, 2 \\
10095 Grugliasco, Torino\\
Italy\\
E-mail: $\mathrm{alessandro.portaluri@unito.it}$

\end{document}